\documentclass[12pt]{amsart}

\numberwithin{equation}{section}
\newcommand{\lra}{\longrightarrow}
\newcommand{\ra}{\to}
\newcommand{\restr}{\mbox{\Large \(|\)\normalsize}}
\def\RR{\mathbf{R}}

\newcommand{\invlim}{\operatorname{\underset{\longleftarrow}{\lim}}}
\newcommand{\dirlim}{\operatorname{\underset{\longrightarrow}{\lim}}}

\def\eps{\epsilon}

\def\lra{\longrightarrow}

\def\ra{\to}

\def\th{\theta}

\def\defeq{:=}

\newcommand{\no}{\noindent}

\def\XXint#1#2#3{{\setbox0=\hbox{$#1{#2#3}{\int}$}
     \vcenter{\hbox{$#2#3$}}\kern-.5\wd0}}

\theoremstyle{plain}
\newtheorem{theorem}[equation]{Theorem}

\newtheorem{proposition}[equation]{Proposition}
\newtheorem{lemma}[equation]{Lemma}
\newtheorem{corollary}[equation]{Corollary}

\theoremstyle{remark}

\theoremstyle{definition}
\newtheorem{definition}[equation]{Definition}

\def\RR{\mathbf{R}}

\newcommand{\N}{\mathbb N}

\newcommand{\R}{\mathbf{R}}

\newcommand{\lla}{\longleftarrow}

\def\ra{\to}

\def\th{\theta}

\def\defeq{:=}

\begin{document}

\begin{abstract}
In this paper we clarify the relation between
inverse systems, the Radon-Nikodym property, the
Asymptotic Norming Property of James-Ho \cite{jamesho},
and the GFDA spaces introduced in \cite{gfda}.  
\end{abstract}

\title[Inverse limits and the RNP]
{Characterization of the Radon-Nikodym Property in terms of inverse
limits}

\date{\today}
\author{Jeff Cheeger}
\address{J.C.\,: Courant Institute of Mathematical Sciences\\
       251 Mercer Street\\
       New York, NY 10012}
\author{Bruce Kleiner}
\address{B.K.\,: Mathematics Department\\
         Yale University\\
             New Haven, CT 06520}
\thanks{The first author was partially supported by NSF Grant
DMS 0105128 and the second by NSF Grant DMS 0701515}
\maketitle

\section{Introduction}

A  Banach space $V$ is said to have the 
{\it Radon-Nikodym Property}  (RNP) if 
every Lipschitz map $f:\RR\to V$ is differentiable almost 
everywhere. 
By now, there are a number of
characterizations of  Banach spaces with the RNP,
 the study of which
 goes back to Gelfand \cite{gelfand}; for additional references
and discussion, see
\cite[Chapter 5]{benlin}, \cite{GhMau1}. 
Of particular interest here is the characterization
of the RNP in terms of the Asymptotic Norming Property; 
\cite{jamesho,GhMau1}.

 In this paper we will show that a variant of the
GFDA property introduced in \cite{gfda} is actually equivalent
to the Asymptotic Norming property of James-Ho, and 
hence by \cite{jamesho,GhMau1}, is equivalent to the
RNP.  In addition, we observe that the GFDA spaces
of \cite{gfda} are just spaces which are isomorphic to 
a separable dual space. 

\begin{definition}
An inverse system 
\begin{equation}
\label{is}
W_1\stackrel{\th_1}{\lla} W_2\stackrel{\th_2}{\lla}\ldots
\stackrel{\th_{i-1}}{\lla} W_{i}\stackrel{\th_{i}}{\lla}\ldots\, ,
\end{equation}
is {\em standard} if the $W_i$'s are finite dimensional Banach 
spaces and the  $\th_i$'s are linear maps of norm $\leq 1$. 
We let $\pi_j:\invlim\, W_i\ra W_j$ denote the projection map.

\end{definition}

\begin{definition}
\label{defdp}
Let $\{(W_i,\th_i)\}$ be
a standard inverse system and 
$V\subset \invlim\, W_i$ be a subspace.
The pair $(\invlim W_i,V)$ has the {\em Determining
Property}  
if a sequence $\{v_k\}\subset V$ converges strongly provided 
the projected sequences $\{\pi_j(v_k)\}\subset W_j$ converge
for every $j$, the sequence
$\{\|v_k\|\}$ is bounded, and the convergence
$\|\pi_j(v_k)\|\to \|v_k\|$ is uniform in $k$.
A Banach space $U$ has the {\em Determining Property}
if there is a pair $(\invlim\,W_i,V)$ with Determining
Property, such that $V$ is isomorphic to $U$.
\end{definition}

 We have:

\begin{theorem}
\label{thmmain}
A separable Banach space has the RNP if and only
it has the Determining Property.
\end{theorem}
Since a Banach space has the 
RNP if and only if every 
separable subspace has the RNP,
Theorem \ref{thmmain}  
yields a characterization of the RNP 
for nonseparable Banach spaces as well. 

To prove the theorem, we first observe
in Proposition \ref{dualem} that the inverse limit
$\invlim\,W_i$ is the dual space of a separable Banach
space.  Then, by a completely elementary argument, we show
that a Banach space  has the 
Determining Property if and only if it has the
Asymptotic Norming Property (ANP) of James-Ho
\cite{jamesho}.  Since
a separable Banach space $U$ has the RNP if and only if
it has the ANP \cite{jamesho,GhMau1}, the theorem
follows.  We remark that there is a simple direct
proof that if $V$ has the ANP (or the Determining Property), then
every Lipschitz map $f:\R\ra V$ is differentiable
almost everywhere,
see \cite{dppi}.

Characterizations of the RNP using inverse limits
are  
useful for applications; 
see \cite{gfda}, the discussion below
concerning metric measure spaces, and \cite{dppi}.

\vskip10mm

\no
{\bf Relation with previous work.}

In slightly different language, 
our earlier paper \cite{gfda} also
considered pairs $(\invlim W_i,V)$, where $\invlim W_i$
is the inverse limit of a standard inverse system,
and $V\subset\invlim W_i$ is a closed subspace.
A  Good
Finite Dimensional Approximation (GFDA)
of a Banach space $V$, a notion introduced in \cite{gfda},
is a pair $(\invlim W_i,V)$ 
with the Determining Property such that 
$\pi_i\restr_V:V\ra W_i$ is a quotient map 
for every $i$.

 It follows immediately from
Lemma 3.8 of \cite{gfda} that if $(\invlim W_i,V)$  
is a GFDA of $V$, then
$V=\invlim\,W_i$.
Since such  inverse limits are dual
spaces by Proposition \ref{dualem}, $V$ is a separable dual space 
in this case.
Conversely,  
using the Kadec-Klee renorming 
Lemma \cite{kadec,klee}, it was shown in \cite{gfda} that 
every separable dual 
space is isomorphic to a Banach space
which admits a GFDA. 
Thus, a Banach space admits a GFDA if and only
if it is isomorphic to a separable dual space.

\bigskip

\noindent
{\bf Applications to metric measure spaces.}

We will call a metric measure space $(X,\mu)$ 
a {\em PI space} if the measure is doubling,
and a Poincar\'e inequality holds in the sense of upper gradients
\cite{heko,cheeger}. 
In \cite{gfda},
differentiation and  bi-Lipschitz non-embedding
theorems were proved for  maps $f:X\to V$ 
from  PI spaces into GFDA
targets $V$, 
generalizing results of \cite{cheeger} for finite 
dimensional targets.
As explained above, 
it turns out that these targets are just separable dual spaces,
up to isomorphism.

As an application of the inverse limit framework
and the equivalence between the ANP and RNP,
we will show in \cite{dppi} that the 
differentiation theorem  \cite[Theorem 4.1]{gfda}
and bi-Lipschitz non-embedding theorem
\cite[Theorem 5.1]{gfda} hold whenever the target has the RNP.

\bigskip
\noindent
{\bf Acknowledgement.}  
We are very grateful to Bill Johnson for sharing an observation 
which helped give rise to this paper. We are much indebted
to Nigel Kalton for immediately
catching a serious error in an earlier version.

\section{Inverse systems}
\label{secinversesystems}
In this section, we recall some basic facts 
concerning direct and inverse systems, and
the duality between them.  Then we show that inverse
limits of standard inverse systems are precisely
duals of separable spaces.

The following conventions will be in force throughout the remainder
of the paper.

\begin{definition}
\label{direct1}
An {\it standard direct system}   is  a sequence
of finite dimensional Banach spaces $\{E_i\}$ and $1$-Lipschitz 
linear maps
$\iota_i:E_i\to E_{i+1}$.
\end{definition}
\vskip1mm

\begin{definition}
\label{inverse1}
An {\it standard inverse system} is a 
sequence of finite dimensional Banach spaces $\{W_i\}$
and $1$-Lipschitz linear maps $\theta_i:W_{i+1}\to W_i$.
\end{definition}
\vskip1mm

\begin{definition}
\label{direct}
A standard direct system is {\it isometrically injective}  
if the maps
$\iota_i:E_i\to E_{i+1}$ are isometric injections.
\end{definition}
\vskip1mm

\begin{definition}
\label{inverse}
A standard inverse system is {\it quotient} if the 
maps $\theta_i:W_{i+1}\to W_i$ are quotient maps.
\end{definition}
\vskip1mm

By a {\it quotient map} of normed spaces, we mean a 
surjective map
$\pi:U\ra V$
for which the norm on the target is the quotient norm, i.e. for 
every $v\in V$, 
$$
\|v\|=\inf\{\;\|u\|\mid u\in\pi^{-1}(v)\}.
$$

We will refer to the maps  $\iota_i$ and $\theta_i$ as {\it bonding 
maps}.

There is a duality between the objects in Definitions \ref{direct1} 
and \ref{inverse1},
respectively,  \ref{direct}
and  \ref{inverse}: if $\{(E_i,\iota_i)\}$ is a standard 
direct system, then $\{(E_i^*,\iota_i^*)\}$ is a standard inverse
system and conversely; similarly, isometrically injective
direct systems are dual to quotient systems.   To see this, one 
uses the facts that 
the adjoint of a $1$-Lipschitz map of Banach spaces
is $1$-Lipschitz and the the adjoint of
an isometric embedding is a quotient map.
(This follows from the Hahn-Banach theorem.)  In particular,
since the spaces in our systems are assumed
to be finite dimensional (hence reflexive)  every  inverse system 
arises as the dual of its dual direct system and conversely. The 
same
holds for quotient inverse systems.

We now recall the definitions of direct and inverse limits.

Given a standard direct system $\{(E_i,\iota_i)\}$ we form the 
direct
limit Banach space $ \dirlim E_i$ as follows.  We begin
with the disjoint union $\sqcup_i\, E_i$, and   
declare two elements $e\in E_i$, $e'\in E_{i'}$ to be equivalent 
if their
images in $E_j$ coincide for some $j\geq \max\{i,i'\}$.
Since the bonding maps are $1$-Lipschitz, the set of 
equivalence classes inherits an obvious vector space structure
with a pseudo-norm.
The direct limit $\dirlim\, E_i$
is defined to be the completion of the quotient of this
space by the closed subspace of elements whose pseudo-norm is zero.
Clearly, there are $1$-Lipschitz maps
$$
\tau_i:E_i\to\dirlim E_i\,,
$$
which in the case of isometrically injective direct systems,
 are isometric injections.  
The union $\bigcup_i \;\tau_i(E_i)$ is dense in $ \dirlim E_i$.

The inverse limit $\invlim W_i$ of a standard inverse
system  $\{(W_i,\theta_i)\}$ is defined as follows.
The underlying set consists of the collection of 
elements $(w_i)\in \prod_i\, W_i$ 
which are compatible with the bonding maps, i.e.
 $\theta_{i}(w_{i})=w_{i-1}$ for all $i$, and
which satisfy $\sup_i \|w_i\| < \infty$.  This  is
equipped with the obvious vector space structure and the norm
\begin{equation}
\label{normdef}
\|\{w_i\}\|\defeq \lim_{j\ra\infty}\|w_j\|\, .
\end{equation}
The  map
\begin{equation}
\label{eqn:pij}
\pi_j:\invlim W_i\to W_j
\end{equation}
given by
$$
\pi_j(\{w_i\})=w_j
$$
is $1$-Lipschitz, and
$$
\lim_{j\to\infty}\|\pi_j(\{w_i\})\|=\|\{w_i\}\|\, .
$$

 An inverse limit $\invlim W_i$ has a natural {\it inverse limit topology}, 
namely the
weakest topology such that every projection map 
$\pi_j:\invlim W_i\to W_j$
is continuous.  Thus a sequence $\{v_k\}\subset \invlim W_i$
converges in the inverse limit topology to $v\in \invlim W_i$ if and only if for every $i$,
 we have
$\pi_i(v_k)\ra \pi_i(v)$
as $k\ra \infty$.

If $\{v_k\}\subset \invlim W_i$ and
  $\{v_k\}\stackrel{invlim}{\longrightarrow}v\in \invlim W_i$, then
\begin{equation}
\label{lowersemicont}
\|v\| \;\leq\; \liminf_k\,\|v_k\|  \, .
\end{equation}
Also, every norm bounded sequence $\{v_k\}\subset\invlim W_i$
 has a   subsequence which converges with respect to the 
inverse limit topology; this follows from a diagonal
argument, because $\{\pi_i(v_k)\}$
is contained in a compact subset of $W_i$, for all $i$.

\bigskip
\begin{proposition}
\label{dualem}
 Given a standard 
inverse system $\{(W_i,\theta_i)\}$, there is an
isometric isomorphism
\begin{equation}
\label{canemb}
C:\invlim\, W_i\equiv (\dirlim W_i^*)^*\, .
\end{equation}
In particular, $\invlim W_i$ is the  dual 
of the separable Banach space $\dirlim W_i^*$.
\end{proposition}
\proof
Pick a compatible sequence $(x_i)\in \invlim W_i$.
We get a map 
$$
\sqcup\; W_j^*\ra \R
$$
by sending $\phi\in W_j^*$ to $\phi(x_j)$; because
$(x_i)$ is compatible with bonding maps and
$$
|\phi(x_j)|\leq \|\phi\|\,\|x_j\|\leq \|\phi\|\,\|\{x_j\}\|,
$$
this defines a linear functional of norm $\leq \|\{x_j\}\|$
on $\dirlim\,W_i^*$.  Therefore we get a $1$-Lipschitz
map
$$
C:\invlim W_i\lra \left(\dirlim W_i^*\right)^*.
$$

We now verify that $C$ is an isometry.

Pick $(x_i)\in \invlim W_i$, and choose
$n\in\N$ such that   $\|x_n\|\geq \|(x_i)\|-\eps$.  
If $\phi\in W_n^*$ has norm $1$ and $\phi(x_n)=\|x_n\|$,
then 
$$
\|C((x_i))\|\,\|\tau_n(\phi)\|\geq C((x_i))(\tau_n(\phi))=\phi(x_n)
=\|x_n\|\geq \|(x_i)\|-\eps,
$$
where $\tau_n:W_n^*\ra\dirlim W_i^*$
is the canonical $1$-Lipschitz map described above.
This shows that $C$ is an isometric embedding.

If $\Phi\in (\dirlim W_i^*)^*$, then we define
$\Phi_i\in W_i^{**}=W_i$ to be the 
composition
$$
W_i^*\lra\dirlim W_i^*\stackrel{\Phi}{\lra}\,\R.
$$
This  defines
a compatible sequence $(\Phi_i)\in\invlim W_i$, such
that $\|(\Phi_i)\|=\|\Phi\|$ and $C((\Phi_i))=\Phi$.
Hence $C$ is onto.
\qed

\begin{corollary}
\label{coryisdirlim}

\mbox{}

\no
{\rm 1)}
A separable Banach space $Y$ is isomorphic to the direct limit of 
an isometrically injective direct system $(E_i,\iota_i)$.

\no
{\rm 2)}
The dual space $Y^*$ of the separable Banach space $Y$ (as in 1))
is isometric to the inverse limit
$\invlim\, E_i^*$ of the 
a quotient inverse system $\{(E_i^*,\iota_i^*)\}$.
\end{corollary}
\proof
To see that 1) holds,  start with  a countable increasing sequence
$E_1\subset E_2\subset\cdots\subset Y$ 
of finite dimensional subspaces whose union is dense
in $Y$,
and take the bonding maps $\iota_i:E_i\ra E_{i+1}$ to be the
inclusions.  Clearly the inclusion maps
$E_i\ra Y$ induce an isometry $\dirlim E_i\ra Y$.

Assertion 2) follows from 1) and Proposition \ref{dualem}.
\qed

\bigskip
\bigskip
Let $C$ be the isometry in Proposition \ref{dualem}.
\begin{lemma}
\label{dualemproperties}

\mbox{}

\no
{\rm 1)}  Suppose $\{v_k\}\subset \invlim\, W_i$ is a sequence such that
 $\{C(v_k)\}\subset (\dirlim W_i^*)^*$  weak* converges to 
some $y\in (\dirlim W_i^*)^*$.   Then  $\{v_k\}$
is   convergent with respect to the inverse limit topology, and its limit 
$v_\infty\in \invlim\, W_i$ 
satisfies $C(v_\infty)=y$; in particular, $y\in C(\invlim\, W_i)$. 

\no
{\rm 2)} If $\{v_k\}\subset \invlim\, W_i$ converges in the
inverse limit topology, and has  uniformly bounded norm, then 
 $\{C(v_k)\}$ is weak* convergent.
\end{lemma}
\proof
Assertions 1) and 2) follow readily from the assumption
that the $W_i$ are finite dimensional together
with the density of compatible sequences  in $\invlim\, W_i$.
\qed

\section{The proof of Theorem \ref{thmmain}}

The proof of Theorem \ref{thmmain} is based on the 
Asymptotic Norming Property, which we now recall.

Let $Y$ denote a separable Banach space and $V\subset Y^*$ 
a separable subspace 
of its  dual.  (Here $Y^*$  need not be separable.)

\medskip
\begin{definition}
\label{def:anp}
The pair $(Y^*,V)$ has the {\em Asymptotic  Norming  
Property} (ANP)  if a sequence $\{v_k\}\subset V$
converges strongly provided it is  weak* convergent
and the sequence of norms $\{\|v_k\|\}$ converges to 
the norm of the weak* limit.

A Banach space $U$ is said to have the 
{\em Asymptotic Norming Property} if there is a 
pair $(Y^*,V)$ with the 
ANP  such that $U$ is isomorphic to $V$.
\end{definition}

\begin{theorem}[\cite{jamesho,GhMau1}]
\label{thmgm}
For separable Banach spaces, the RNP
is equivalent to  the ANP.
\end{theorem}

Hence  to prove
Theorem \ref{thmmain}, it suffices to show that
for separable Banach spaces, the ANP is equivalent
to the Determining Property.  By Corollary
\ref{coryisdirlim}, every separable Banach space $Y$
is isometric to 
the direct limit of a standard direct system, and
$Y^*$ is isometric to the inverse limit of the dual
inverse system.  Hence the
proof of Theorem \ref{thmmain} reduces to:

\begin{proposition}
\label{propanpequalsdp}
Let $\{(W_i,\th_i)\}$ be a standard inverse system,
and $V$ be a closed separable subspace of $\invlim W_i$.
Then the pair $(\invlim W_i,V)$ has the ANP if and only
if it has the Determining Property.  Here we are identifying
$\invlim W_i$ with the dual of $\dirlim W_i^*$, see 
Proposition \ref{dualem}.
\end{proposition}
\proof
Let $\{v_k\}\subset V$ be a sequence with bounded
norm.  By Lemma \ref{dualemproperties}, the sequence 
$\{v_k\}$ is weak* convergent if and only if it
converges  in the inverse limit topology.  
Therefore, to prove the equivalence
of the ANP and the Determining Property for the
pair $(\invlim W_i,V)$, it suffices to show that
when 
\begin{equation}
\label{anp1}
v_k\stackrel{w^*}{\longrightarrow} w\in\invlim \, W_i\, ,
\end{equation}
the sequence of norms $\{\|v_k\|\}$ converges to the
$\|w\|$ if and only if the 
convergence $\|\pi_j(v_k)\|\ra \|v_k\|$
is uniform in $k$.
Although
this is completely elementary, we will write out the details.

We have
\begin{equation}
\label{3}
\|v_k\|-\|w\|=(\|v_k\|-\|\pi_i(v_k)\|)
+(\|\pi_i(v_k)\|-\|\pi_i(w)\|)+(\|\pi_i(w)\|-\|w\|)\, .
\end{equation}

Assume first that $\lim_{k\to\infty}\|v_k\|=\|w\|$.  
Given $\epsilon>0$, there
exists $I_1$ such that $\|w\|-\|\pi_i(w)\|<\epsilon/3$, 
for $i\geq I_1$.  By (\ref{anp1})
there exists $K_1$ such that 
$\|\pi_{I_1}(v_k)-\pi_{I_1}(w)\|<\epsilon/3$, for
$k\geq K_1$.  Also, there exists $K_2$ such that 
$\big | \|v_k\|-\|w\|\big |<\epsilon/3$,
if $k\geq K_2$.  Set $K=\max(K_1,K_2)$. 

 From (\ref{3}), with $i=I_1$, we get 
$\|v_k\|-\|\pi_{I_1}(v_k)\|<\epsilon$, for all $k\geq K$.
Since, $\|v_k\|-\|\pi_i(v_k)\|$ is a nonnegative decreasing 
function of $i$,
this implies, $\|v_k\|-\|\pi_{i}(v_k)\|<\epsilon$, for all
 $i\geq I_1$, $k\geq K$.

 Finally, there exists $I_2$ such that 
$\|v_k\|-\|\pi_i(v_k)\|<\eps$ 
for all $i\geq I_2$,
$k=1,\ldots, K-1$, Thus, if $i\geq \max(I_1,I_2)$ then
$\|v_k\|-\|\pi_i(v_k)\|<\epsilon$, for all $k$.

Conversely, suppose the convergence $\|\pi_i(v_k)\|\to \|v_k\|$ 
is uniform in $k$.
Given $\epsilon>0$, there exists $I$ such that 
$\|v_k\|-\|\pi_i(v_k)\|<\epsilon/3$, for 
$i\geq I$ and all $k$. Also, there exists $I_1$ such that
 $\|w\|-\|\pi_i(w)\|<\epsilon/3$, 
for $i\geq I_1$. Set $I'=\max(I,I_1)$. By (\ref{anp1}), 
there exists $K$ such that
$\|\pi_{I'}(v_k)-\pi_{I'}(w)\|<\epsilon/3$.

 From (\ref{3}), with $i=I'$, we get 
$\big | \|v_k|-\|w\|\big |<\epsilon$, for all $k\geq K$.
\qed

\section{A variant of the Determining Property}

In this section we discuss a variant of the
Determining Property, which was introduced in 
\cite{gfda} (with a different name).  
A compactness argument implies
that it is equivalent
to Definition \ref{defdp}, see Proposition 
\ref{uniformityclass}.

For the remainder of this section, we fix
a standard inverse system $\{(W_i,\th_i)\}$
and a closed subspace $V\subset\invlim\, W_i$.

\begin{definition}
\label{defepsdetermining}
A  positive nonincreasing 
finite sequence $1\geq \rho_1\geq\ldots \geq\rho_N$ 
is {\em $\eps$-determining} if for any pair 
$v,v'\in V$, the conditions
\begin{equation}
\label{e:epsdet.1}
{}
\end{equation}
$$
\|v\|-\|\pi_i(v)\|<\rho_i\cdot  \|v\|, 
\qquad  \|v'\|-\|\pi_i(v')\|
<\rho_i\cdot \|v'\|, \qquad 1\leq i\leq N\, ,
$$
and 
\begin{equation}
\label{e:epsdet.2}
\|\pi_N(v)-\pi_N(v')\|<N^{-1}\cdot\max(\|v\|,\|v'\|)\, ,
\end{equation}
imply
\begin{equation}
\label{e:epsdet.3}
\|v-v'\|<\eps\cdot\max(\|v\|,\|v'\|)\, .
\end{equation}
\end{definition}
\no
Observe that by dividing by $\max(\|v\|,\|v'\|)$, 
it suffices to consider
pairs $v,v'$ for which $\max(\|v\|,\|v'\|)=1$.

This leads to the alternate definition of the Determining
Property:
\begin{definition}
\label{defdp2}
The pair $(\invlim W_i,V)$ 
has the  {\it Determining Property}
if for every $\eps>0$ and
 every infinite nonincreasing sequence
$$
1\geq \rho_1\geq \ldots\geq\rho_i\geq\ldots
$$
with $\rho_i\ra 0$,  
some finite initial segment 
$\rho_1\geq\ldots\geq\rho_N$
is $\eps$-determining.
\end{definition}

\begin{proposition} 
\label{uniformityclass}
\hspace{.2in}The pair
$(\invlim W_i,V)$ 
satisfies Definition \ref{defdp} if and only if
it satisfies Definition \ref{defdp2}.
\end{proposition}
\begin{proof}
 First we show  that the property in Definition \ref{defdp2} 
implies the property in
Definition \ref{defdp}.
 So assume that the sequence
$\{\|v_k\|\}$ is bounded and the convergence, 
$\|\pi_i(v_k)\|\to \|v_k\|$ is uniform in
$k$. 

Suppose that there exists a sequence, 
a positive sequence, $\rho_i\searrow 0$, such that
$\|v_k\|-\|\pi_i(v_k)\|\leq \rho_i$. By applying the condition 
in Definition \ref{defdp2} to this sequence and using 
convergence in the inverse limit topology
together with (\ref{e:epsdet.2})
it is clear from (\ref{e:epsdet.3}) that we obtain strong 
convergence.

Without loss of essential loss of generality, we can assume
$\|v_k\|\leq 1$ for all $k$. Since the convergence, 
$\|\pi_i(v_k)\|\to \|v_k\|$ is uniform in $k$, it follows that
there exists a strictly increasing sequence,
 \hbox{$N_1<N_2<\ldots$},
such that for all $k$, we have
$$
\|v_k\|-\|\pi_{N_\ell}(v_k)\|<\frac{1}{\ell}\, .
$$
Then $\|v_k\|-\|\pi_i(v_k)\|\leq \rho_i$, for the sequence, 
$\rho_i$ given
by 
$$
\rho_i=\frac{1}{\ell}\qquad (N_\ell\leq i< N_{\ell+1})\, .
$$

Conversely, suppose that the property in Definition \ref{defdp} 
holds,
but not the property in Definition \ref{defdp2}.
Then for
some decreasing sequence
 $\{\rho_i\}\subset (0,\infty)$ with $\rho_i\ra 0$,
and some $\eps>0$, there are sequences $\{v_k\},\{v_k'\}\subset V$, 
such that for all $k<\infty$,
\begin{equation}
\label{norm}
\|v_k\|,\, \|v'_k\|\leq 1\, ,
\end{equation}
\begin{equation}
\label{unif}
\max\left(\|v_k\|-\|\pi_i(v_k)\|,\,\;\|v'_k\|
-\|\pi_i(v'_k)\|\right)<\rho_i\;\;
\mbox{for}\; 1\leq i\leq k\, ,
\end{equation}
\begin{equation}
\label{close}
\|\pi_i(v_k)-\pi_i(v'_k)\|<\frac{1}{k}\, ,
\end{equation}
\begin{equation}
\label{sep}\|v_k-v_k'\|\geq \eps\, .
\end{equation}

 By the Banach-Alaoglu theorem, we can
pass to weak${}^*$ convergent subsequences, with 
respective limits
$v_\infty$ and $v_\infty '$. From (\ref{close}), it follows that
$v_\infty=v_\infty'$.

It follows from 
(\ref{norm}), (\ref{unif}), that the sequences,  
$\|v_k\|$, $\|v_k'\|$,
are bounded and the convergence $\|\pi_i(v_k)\|\to \|v_k\|$,
$\|\pi_i(v_k')\|\to \|v_k'\|$ is uniform in $k$. 
Since we assume the 
property in Definition \ref{defdp}, it follows
$v_k\to v_\infty$, $v_k'\to v'_\infty$,  is actually strong.
Since, $v_\infty=v'_\infty$, this contradicts (\ref{sep}).
\end{proof}

We remark that proof of the  implication Definition 
\ref{defdp}$\implies$ Definition \ref{defdp2}
is similar to the proof of Proposition 3.11 in
\cite{gfda}.

\section{GFDA versus ANP}
\label{secearlierresults}

We conclude with some remarks about the relation between
the ANP and GFDA's.

Suppose $Y$ is a separable Banach space and $(Y^*,V)$
has the ANP.  By Lemma \ref{coryisdirlim}, we may realize
$Y^*$ -- up to isometry -- as the inverse limit of 
a quotient system $\{(W_i,\th_i)\}$.  

Viewing $V$ as
a subspace of $\invlim W_i$, one might be tempted to 
modify the inverse system to produce a GFDA of $V$.
For instance, one could 
restrict the projection maps
$\pi_j:\invlim W_i\ra W_j$
to $V$, and replace $W_j$ with $\pi_j(V)\subset W_j$.
However, the resulting maps $\pi_j\restr_V:V\ra \pi_j(V)$
will usually not be quotient
maps.  One could also try renorming the spaces $\pi_j(V)\subset
W_j$ so that the restrictions $\pi_j\restr_V:V\ra\pi_j(V)$
become quotient maps.  This will typically destroy the
Determining Property, however.
In any case,  $V$ will not admit any GFDA unless
it is a separable dual space, whereas many Banach
spaces with the RNP are not separable
dual spaces.

\bibliography{gfda}
\bibliographystyle{alpha}
\addcontentsline{toc}{subsection}{References}

\end{document}